\documentclass[11pt]{article}
\usepackage{amsmath,amssymb,amsthm,amsfonts,amstext,amsbsy,amscd}
\usepackage{mathrsfs}
\usepackage{mathabx}

\usepackage{float}

\usepackage{bm}
\usepackage{subcaption}
\usepackage{multirow}
\usepackage{caption}
\usepackage{tabularx}
\usepackage{booktabs}
\usepackage{multirow}
\usepackage{graphicx}
\usepackage{multicol}
\usepackage[ruled,vlined]{algorithm2e}
\usepackage{comment}
\usepackage{enumerate}
\usepackage{latexsym}
\usepackage[top=3cm, bottom=3.5cm, left=3cm, right=3cm]{geometry}

\usepackage[utf8]{inputenc}
\usepackage{dsfont}
\usepackage{color}
\usepackage{tikz}
\RequirePackage[colorlinks,citecolor=blue,urlcolor=blue]{hyperref}

\newcommand{\ind}{\mathbf{1}}
\newtheorem{theorem}{Theorem}

\newtheorem{remark}{Remark}

\newtheorem{proposition}{Proposition}
\renewcommand{\tilde}{\widetilde}

\newcommand{\E}{\mathbb{E}}
\newcommand{\C}{\mathbb{C}}

\newcommand{\R}{\mathbb{R}}
\newcommand{\N}{\mathbb{N}}
\newcommand{\Pro}{\mathbb{P}}

\begin{document}
%%-----------------------------
%%      the top matter
%%-----------------------------
%\title{Stable and Classical Temperated stable Lévy process: properties and trajectory simulation.}
\title{Stable and tempered stable distributions and processes: an overview toward trajectory simulation}
\author{Jalal Taher \footnote{UMR 8100 - Laboratoire de Mathématiques de Versailles, UVSQ, CNRS, Université Paris-Saclay, 78035
Versailles Cedex, France.}}
\date{}

\maketitle
\begin{abstract} 
Stable distributions are a celebrated class of probability laws used in various fields. The $\alpha$-stable process, and its exponentially tempered counterpart, the Classical Tempered Stable (CTS) process, are also prominent examples of Lévy processes. Simulating these processes is critical for many applications, yet it remains computationally challenging, due to their infinite jump activity. This survey provides an overview of the key properties of these objects offering a roadmap for practitioners. The first part is a review of the stability property, sampling algorithms are provided along with numerical illustrations. Then CTS processes are presented, with the Baeumer-Meerschaert algorithm \cite{Baeumer2010TemperedSL} for increment simulation, and a computational analysis is provided with numerical illustrations across different time scales.
\end{abstract}
\section*{Introduction}

Among the many examples of Lévy processes, the $\alpha$-stable process stands out as a central example of processes 
%characterized by \emph{càdlàg} trajectories 
with an infinite jump activity. % Due to their stationary and independent increments,%
Understanding these processes reduces to the study of stable distributions, a topic extensively covered in the literature. Stability was first characterized in 1925 by Levy in \cite{lévy1925calcul} while exploring the properties of the sum of i.i.d. random variables. Gnedenko and Kolmogorov later emphasized the potential applications of stable distributions in \cite{gnedenko_limit_1954}, stating that it deserves "\emph{the most serious attention. It is probable that the scope of applied problems in which they play an essential role will become in due course rather wide}".

There are many reasons why stable distribution and processes are attractive for applications. They encompass a wide range of behaviors; e.g. they include the Cauchy and the Gaussian distributions as special cases. They also arise in the generalized Central Limit Theorem as the only limit of normalized sums of i.i.d. random variables. Additionally, their heavy tails align with studies suggesting the importance of extreme events, providing an alternative to Gaussian models. Therefore, stable distributions and processes have found applications across various fields, including behavior studies \cite{viswanathan_optimizing_1999,reynolds_free-flight_2007,brockmann_scaling_2006}, signal processing \cite{Nolan_filter}, computer science \cite{Terdik_2009}. It also holds a particular place in finance since Mandelbrot in \cite{mandelbrot_variation_1997} discarded the historical Gaussian assumptions for asset returns and inaugurated the use of stable distributions. It has been central in pricing options \cite{kaplan2012frontiers} and modeling the price of commodities (e.g. electricity \cite{weron_heavy_2005}, agricultural \cite{jin_heavy-tailed_2007}). There also has been a renewed interest in heavy-tail distributions in more recent topics in deep-learning \cite{gurbuzbalaban_heavy-tail_2021,jung_-stable_2023}. For further details on applications, refer to \cite{nolan_univariate_2021,zolotarev_one-dimensional_1986,cont_financial_2003,rachev2011financial}.

However, the infinite variance of non-Gaussian stable distributions (see \eqref{eq:moment_stable}) and the inadequacy of this assumption in various models (see \cite{rachev2011financial}), led to considering alternative distributions.  Methods to circumvent this issue have historically originated in physics, trying to maintain a stable-like distribution in some central region while lightening the tails so the variance is finite. Mantegna and Stanley \cite{mantegna_stochastic_1994} introduced the truncated Lévy flight, cutting the density function $f$ of a stable distribution; $f_T(x) \propto f(x) \ind_{|x|\leq T}.$ While this modification yields finite moments, and resembles stable distribution for large $T$, it looses the infinite divisibility property (see Remark \ref{rem:infinite_divisibility}) which is a central tool for computations. Koponen in \cite{koponen_analytic_1995} solved this issue with a smooth exponential cutoff resulting in the Classical Tempered Stable (CTS) distribution. They have been discovered multiple times resulting in various terminology; smoothly truncated Lévy flight (STLF) \cite{koponen_analytic_1995,grabchak_tempered_2016}, Tempered stable (TS) \cite{cont_financial_2003,kuchler_tempered_2013,nolan_univariate_2021}, KoBol in \cite{kobol}, CGMY (for the symmetric version) \cite{carr_fine_2002}, exponentially tilted stable distributions in the sampling literature \cite{hofert_sampling_2011,Devroy_2009}. Notable limit cases include the bilateral Gamma \cite{kuchler_Tappe_2008} and the Variance Gamma \cite{madan_variance_1990} distributions. In this article, we adopt the term Classical Tempered Stable (CTS) following \cite{rachev2011financial}. Finally, in \cite{rosinski_generalized_2010}, Rosinski generalized the CTS distribution with different tempering functions and in higher dimensions; these are referred to as Tempered Stable distribution. The CTS discussed in this article can be seen as a specific subclass in dimension $d=1$ with exponential tempering function (see Chapter 6 of \cite{grabchak_tempered_2016}). The CTS distributions and processes are also widely applied in fields including actuarial science \cite{griffin2013finite,hainaut2008mortality} cellular biology \cite{palmer2008modelling}, computer science \cite{Cao2014}. But it has been extended and considerably developed in financial applications (see \cite{grabchak_tempered_2016, rachev2011financial} for an overview). 

Simulating trajectories of stable and CTS processes is crucial for their use in applications. These processes exemplify infinite activity Lévy processes and are used in many inference models \cite{zhao2009nonparametric,mies2019nonparametric}. They also serve as numerical validation proxies of statistical methods \cite{FigueroLopezNonParametric2024,CeciliaDriftBurst2022,duval2024nonparametric}. In practice, once the simulation of stable and CTS distribution is achieved, a trajectory approximation of the corresponding process can be derived thanks to the independent and stationary increments property. This incorporates a temporal aspect to the problem through the sampling rate $\Delta>0$ which can impact the computational efficiency of algorithms.
For the stable distributions, direct sampling algorithms relating to the work of Chambers, Mallows, Collin, and Stuck \cite{chambers1976method} are available and commonly used in practice \cite{cont_financial_2003,nolan_univariate_2021}. Thanks to the selfsimilarity of $\alpha$-stable process (see \ref{sec:stable_levy_process}), the sampling of $X_\Delta$ can be reduced to that of $X_1$, hence the computational cost of trajectory simulation remains consistent across time scales: it is as costly in high (large $\Delta$) of low (small $\Delta$) frequency (see \ref{sec:simulation_stable}). 

For CTS distribution, simulation is more intricate. Lévy process methods include Compound Poisson approximation complemented with Gaussian approximation of the small jumps (see Section 6.3 of \cite{cont_financial_2003}, \cite{asmussen_approximations_2001,Carpentier2021} for theoretical justifications, Section 4 of \cite{kawai_simulation_2011} for computational details), shot noise series representation (see \cite{Rosiński2001,kawai_simulation_2011} and Section 8.3.2 of \cite{rachev2011financial}), subordination in the symmetrical case (see Proposition 4.1 of \cite{cont_financial_2003} and Section 8.5.1 of \cite{rachev2011financial}). We refer to the textbook by Asmussen and Glynn \cite{asmussen2007stochastic} for a general overview of Lévy process simulation. Alternatively, rejection-based sampling algorithms can be derived from specific properties of CTS \cite{DEVROYE1981547,Baeumer2010TemperedSL,devroye2009random,hofert_sampling_2011}. Kawai and Masuda \cite{kawai_simulation_2011} provide a comprehensive comparison of these methods, highlighting the Baeumer-Meerschaert algorithm \cite{Baeumer2010TemperedSL} as the most efficient and practical to tune for small time increments. It is widely used for its simplicity in recent work involving CTS distributions as in \cite{gajda2010fractional,CHRISTENSEN2022461}.

This article aims to provide hands-on algorithms for simulating increments and trajectories of both $\alpha$-stable and CTS Lévy processes. It also serves as a road map for practitioners to navigate the interlaced concepts and algorithms provided for the CTS processes and its avatars (CGMY, Kobol, Smoothed Levy flight...). It is structured as follows. The first part is a review of stable distributions, presenting the equivalent definitions of stability (Theorem \ref{th:stable_equivalence_def}). Selected properties are outlined with a view toward simulation. Finally, sampling algorithms and numerical illustration are provided in section \ref{sec:simulation_stable}. The second part focuses on the CTS process detailing the Baeumer-Meerschaert algorithm in the bilateral case and discussing computational aspects with illustrative examples.

\textbf{Notations:} For $X$ a random variable, we denote its characteristic function by $\phi_X$, where $\phi_X(u) = \E(e^{iuX})$ for $u \in \R$. If it law  $\Pro_X$ is absolutely continuous with respect to the Lebesgue measure, we denote its density function by $f_X$. The notation $X \stackrel{d}{=} Y$ signifies that random variables $X$ and $Y$ have the same distribution. The upper incomplete Gamma function is defined for $x\geq 0, s>0  $ by  $\Gamma(s,x)= \int_{x}^\infty e^{-x}x^{s-1}dx$ (when $x>0$ it can be extended to an entire function).

\section{Stable distributions and $\alpha$-stable Lévy processes: an overview}
%\textbf{References:}\cite{samoradnitsky_stable_2017,cont_financial_2003,feller_introduction_1991,gnedenko_limit_1954,zolotarev_one-dimensional_1986,nolan_univariate_2021}
\subsection{The stability property}

Four equivalent definitions of stability are presented, each shedding light on a distinct property of stable distributions. The term 'stable' initially refers to \emph{sum stable} as random variables that keep the same shape under summation. Linearly combining two i.i.d. stable variables corresponds to the same distribution as an affine transformation of one of them \eqref{eq:def_stable_2_sum} and so is an $n-$sum \eqref{eq:def_stable_n_sum}. An equivalent approach to stable distributions comes through normalized sums. They constitute the unique limits of normalized sums for i.i.d.  random variables $(X_i)_i$; this corresponds to the Generalized Central Limit Theorem, as outlined in \eqref{eq:def_stable_GCLT}. If $X_i$ has a finite variance, the limit distribution is Gaussian and the classical Central Limit Theorem is recovered. If the variance is infinite, the only limits are in the class of stable distributions, giving a perspective on stability as a generalization of Gaussianity. Finally, stable distributions are within the rich class of infinitely divisible random variables (see \cite{sato_levy_1999} for a detailed presentation on the topic), leading to closed formulas of their characteristic function in  \eqref{eq:def_stable_cf}.
All equivalent definitions are summarized in the following Theorem. 
  
\begin{theorem}\label{th:stable_equivalence_def}
A random variable $X$ is said to have a stable distribution if one of the equivalent following properties is met.
\begin{itemize}
\item \textbf{Sum stability for linear combination.} For any $a,b>0$ there exist two constants $c>0$ and $d\in \R$ such that 
    \begin{align}\label{eq:def_stable_2_sum}
        a X_1 + b X_2 \stackrel{d}{=} cX + d,
    \end{align}
    where $X_1$ and $X_2$ are independent copies of $X$.
    
\item \textbf{Sum stability for $n$ terms sum.} For any $n\in \N$ such that $n\geq 2$, there exist two constants $c_n>0$ and $d_n \in \R$ such that 
    \begin{align}\label{eq:def_stable_n_sum}
        X_1 + X_2 + \cdot + X_n \stackrel{d}{=} c_n X + d_n,
    \end{align}
    where $X_1,X_2, \cdots, X_n$ are independent copies of $X$ and for all $n\in N$ there exists $\alpha \in (0,2]$ called the \emph{stability index} such that $c_n = n^{1/\alpha}$.

\item \textbf{General Central Limit Theorem.} There exist $Y_1,\cdots, Y_n$ i.i.d. random variables and two sequences $(a_n)_{n\in\N}$, $(b_n)_{n\in \N}$, where $\forall n\in \N, a_n>0$ and $b_n \in \R$ such that
\begin{align}\label{eq:def_stable_GCLT}
 a_n\left(\sum_{i=1}^n Y_i\right) - b_n \stackrel{d}{\rightarrow} X,
\end{align}
when $n\rightarrow +\infty$

\item \textbf{Stable characteristic function.} There exist $0<\alpha\leq 2, \sigma \geq 0, \beta \in [-1,1]$ and $\delta \in \R$ such that the characteristic function $\phi_X$ of $X$ is given for all $u\in \R$ by
\begin{align}\label{eq:def_stable_cf}
\phi_X(u) = \begin{cases}
    \exp \left( - \sigma^\alpha |u|^\alpha \left(1- i\beta sign(u) \tan \left(\frac{\pi\alpha}{2} \right) \right) + i\delta u \right)& \text{if } \alpha \neq 1 \\
    \exp \left(- \sigma |u|(1+ i\beta \frac{2}{\pi} sign(u) \log(|u|)) + i \delta u \right) &\text{if } \alpha =1.
\end{cases}
\end{align}
\end{itemize}
\end{theorem}

\begin{figure}[H]
\begin{center}
\includegraphics[scale=0.5]{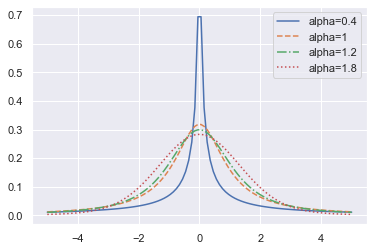}
\includegraphics[scale=0.5]{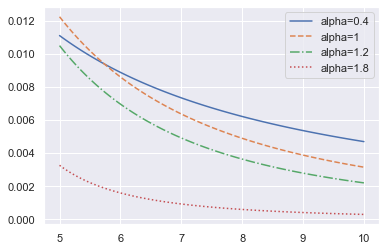}
\caption{Display of the density $f_X(x)$ where $X\sim S_\alpha(1,0,0)$ with varying $\alpha \in \{0.4,1,1.2,1.8\}.$ The index $\alpha$ affects the tail; the smaller $\alpha$ is the heavier is the tail. }
\label{fig:stable_density_alpha}
\end{center}
\end{figure}

\begin{figure}[H]
\begin{center}
\includegraphics[scale=0.5]{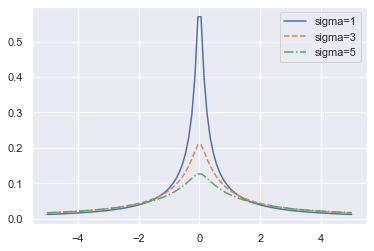}
\includegraphics[scale=0.5]{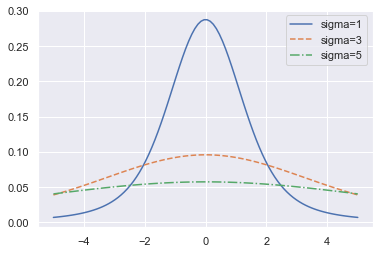}
\caption{Display of the density $f_X(x)$ where $X\sim S_\alpha(1,0,0)$ with varying $\sigma \in \{1,3,5\}$, $\alpha=0.5$ (left) and $\alpha=1.5$ (right). The parameter $\sigma$ has a rescaling effect on the density }
\label{fig:stable_density_sigma}
\end{center}
\end{figure}

\begin{figure}[H]
\begin{center}
\includegraphics[scale=0.5]{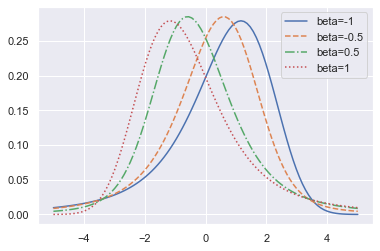}
\includegraphics[scale=0.5]{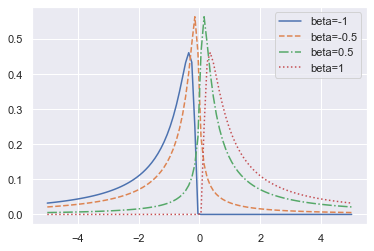}
\caption{Display of the density $f_X(x)$ where $X\sim S_\alpha(1,0,0)$ with varying $\beta \in \{1,3,5\}$, $\alpha=0.5$ (left) and $\alpha=1.5$ (right). The parameter $\beta$ affects the shape of the density. }
\label{fig:stable_density_beta}
\end{center}
\end{figure}

\begin{figure}[H]
\begin{center}
\includegraphics[scale=0.5]{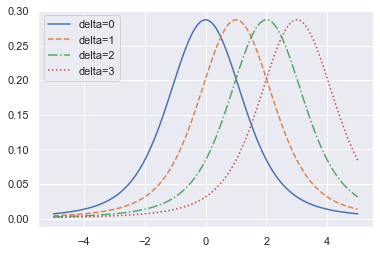}
\caption{Display of the density $f_X(x)$ where $X\sim S_\alpha(1,0,0)$ with varying $\delta \in \{0,1,2,3\}$, $\alpha=0.5$. The parameter $\delta$ has a shifting effect on the density. }
\label{fig:stable_density_delta}
\end{center}
\end{figure}

\begin{remark}
The closed formula \eqref{eq:def_stable_cf} yields the fact that stable distributions are a parametric class of distributions. They are characterized by four parameters.
\begin{itemize}
\item \emph{The stability index} $\alpha \in (0,2]$; it affects the shape of the distribution and its tails (see Figure \ref{fig:stable_density_alpha}), e.g. $\alpha=1$ corresponds to a Cauchy distribution, $\alpha=2$ to a Gaussian distribution 
\item \emph{The skewness parameter} $\beta \in [-1,1]$; the distribution is said to be \emph{totally positively (negatively) skewed} if $\beta=1$ ($\beta=-1$) it also affects the shape of the distribution (see Figure \ref{fig:stable_density_beta}).
\item \emph{The scale parameter} $\sigma>0$; which is not the standard deviation of non-Gaussian stable distributions as the variance is infinite when $\alpha \in (0,2)$ (see \eqref{eq:moment_stable}) (see figure \ref{fig:stable_density_sigma}).
\item \emph{The location parameter} $\delta>0$; it is not the mean (see \eqref{eq:moment_stable}) but has a drifting effect on the distribution (see Figure \ref{fig:stable_density_delta}). 
\end{itemize}
\end{remark}

\textbf{Notation:} In the following a stable distribution will be denoted by $S_\alpha(\sigma,\beta,\delta)$ which corresponds to the standard notation used in \cite{samoradnitsky_stable_2017}. In the literature, several other parametrizations have been used, which might create confusion. This issue is discussed in \cite{HallComedyerror} and described as a 'comedy of errors'. The notation used here corresponds to the $1-$parametrization in \cite{nolan_univariate_2021} and to the parametrization $A$ in \cite{zolotarev_one-dimensional_1986}. See Table in \cite{nolan_univariate_2021} for an exhaustive presentation of the different parametrizations of stable distributions. 

\begin{remark}\label{rem:infinite_divisibility}(\emph{Infinite divisibility}).
A random variable $X$ is infinitely divisible if for all $n\geq 2$ there exist $Y_1^{(n)},\cdots Y_{n}^{(n)}$ independent random variables such that $X \stackrel{d}{=}Y_1^{(n)}+\cdots +Y_{n}^{(n)}$ (e.g. Definition 3.2 in \cite{cont_financial_2003}). Thanks to the Lévy-Khintchine Theorem (Theorem 3.2 in \cite{cont_financial_2003}), the characteristic function of infinitely divisible random variables can be expressed using generic integral terms. They rely on the Lévy triplet $(b,A,\nu)$ where $A\geq 0$, $c\in \R$ and $\nu$ is a $\sigma-$finite measure satisfying $\int_{\R} (x^2\wedge1) \nu(dx) <\infty$. For $u \in \R$ the characteristic function of any infinitely divisible random variable is given by
\begin{align}\label{eq:levy_khintchine}
\phi_X(u) = \exp \left(iuc - \frac{\sigma^2}{2}u^2 + \int_{\R}(e^{iux}-1-iux\ind_{[-1,1]}) \nu(dx) \right).
\end{align}
 Let $X$ be a stable distribution. It follows from \eqref{eq:def_stable_n_sum} that for all $n\geq 2$ there exist $c_n>0$ and $d_n \in \R$ such that $X_1 + \cdots + X_n \stackrel{d}{=} c_n X + d_n.$
Choosing $Y_i^{(n)}= \frac{1}{c_n}\left(X_i - \frac{d_n}{n} \right)$ for $i \in \{1,\cdots,n\},$ ensures that $X$ is infinitely divisible. For stable distributions, the Lévy measure is absolutely continuous with respect to the Lebesgue measure and takes the following form (see \cite{gnedenko_limit_1954,sato_levy_1999}) 
 \begin{align}\label{eq:stable_levymeasure}
        \frac{\nu(dx)}{dx} = \frac{P}{x^{1+\alpha}} \ind_{x>0} + \frac{Q}{|x|^{1+\alpha}} \ind_{x<0},
    \end{align}
where $P,Q\geq 0$ and $P+Q>0$. If $\alpha=2$ the Lévy measure $\nu$ is null and $X$ is a rescaled Gaussian random variable.
\end{remark}

\begin{proof}[Sketch of the proof of Theorem \ref{th:stable_equivalence_def}]
One can show that  \eqref{eq:def_stable_2_sum},\eqref{eq:def_stable_GCLT} and \eqref{eq:def_stable_cf} are respectively equivalent to \eqref{eq:def_stable_n_sum}. For further details, we refer to \cite{feller_introduction_1991,gnedenko_limit_1954,nolan_univariate_2021}.
\begin{itemize}
\item (\eqref{eq:def_stable_2_sum} $\iff$ \eqref{eq:def_stable_n_sum}). The direct implication is achieved by induction. The converse (\eqref{eq:def_stable_n_sum} $\implies$ \eqref{eq:def_stable_2_sum}) leverages the fact that $c_n = n^{1/\alpha}$ (see Theorem 1 in VI.1.1 of \cite{feller_introduction_1991}). For simplicity, we assume that $X$ is strictly stable (i.e. $d_n=0$). If $S_n = \sum_{i=1}^n X_i$ then by regrouping the terms we derive that $S_{n+m} = S_n + (S_{n+m} - S_n). $
Since $S_{n+m} = c_{m+n}X$, by independence we get that $$(n+m)^{1/\alpha} X\stackrel{d}{=} n^{1/\alpha} X_1 + m^{1/\alpha} X_2 .$$
This can be extended by replacing $n$ and $m$, with rational numbers and subsequently with real numbers using a continuity argument. (see \cite{feller_introduction_1991} for more details)

\item (\eqref{eq:def_stable_GCLT} $\iff$ \eqref{eq:def_stable_n_sum}). The implication \eqref{eq:def_stable_n_sum} $\implies $ \eqref{eq:def_stable_GCLT} is straightforward taking $a_n = n^{1/\alpha}$ and $b_n = \frac{d_n}{a_n}$. The conserve follows as a consequence of the Convergence of Type Theorem (see Theorem 3.11 in \cite{nolan_univariate_2021}).

\item (\eqref{eq:def_stable_cf} $\iff$ \eqref{eq:def_stable_n_sum}).  Assuming that $X$ satisfies \eqref{eq:def_stable_cf}, then it can be shown that \eqref{eq:def_stable_n_sum} holds with $c_n=n^{1/\alpha}$ and $d_n = (n-n^{1/\alpha}) \delta$. The converse relies on the infinite divisibility of $X$. The closed formulas for characteristic function are achieved using the Lévy Khintchine formula \eqref{eq:levy_khintchine}, the Lévy measure expression \eqref{eq:stable_levymeasure} along with the complex integrals in Lemma 14.11 of \cite{sato_levy_1999}.
\end{itemize}
\end{proof}

\subsubsection*{General distributional properties}
In this section, we present some important properties of stable distributions. First, the $S_\alpha(\sigma,\beta,\delta)$ parametrization is endowed with simple formulas when rescaling and shifting allowing major simplifications when it comes to simulations (Proposition \ref{prop:sum_CL_stable}). Moreover, asymptotic estimates of the tail behavior of stable distribution  (Proposition \ref{th:tail_stable}) emphasize their utility as proxies for heavy-tailed distributions, while also inducing infinite variance as a counterpart \eqref{eq:moment_stable}. An extensive presentation of the many other properties of stable distributions can be found in the following monographs \cite{samoradnitsky_stable_2017,nolan_univariate_2021,sato_levy_1999,zolotarev_one-dimensional_1986} among others.

The parametric form of the characteristic function \eqref{eq:def_stable_cf}, implies that both rescaling and shifting preserve the stability property. Moreover, under the same stability index assumption, linear combinations of independent stable random variables also remain stable.

\begin{proposition}[Proposition 1.4 of \cite{nolan_univariate_2021}]\label{prop:sum_CL_stable}
Let $X$ and $Y$ be independent stable random variables such that $X\sim S_\alpha (\sigma_0,\beta_0,\delta_0)$ and $Y\sim S_\alpha(\sigma_1,\beta_1,\delta_1)$ where $\alpha \in (0,2]$, $\beta_i \in [-1,1]$ and $\delta_i \in \R$ and $i\in \{1,2\}.$
\begin{enumerate}
    \item For any $a\neq 0, b\in \R$ we have that
    \begin{align*}
        aX + b \sim \begin{cases}
            S_\alpha \left(|a| \sigma_0, sign(a)\beta_0, a\delta_0 +b\right) &\text{ if } \alpha \neq 1\\
            S_1\left( |a| \sigma_0, sign(a)\beta_0, a \delta_0 + b- \frac{2}{\pi}\beta_0 a \log(|a|)\right) &\text{ if } \alpha=1.
        \end{cases}
    \end{align*}

    \item The sum $X+Y \sim S_\alpha \left(\sigma_1^\alpha + \sigma_2^\alpha, \frac{\beta_0 \sigma_0^\alpha + \beta_1 \sigma_1^\alpha}{\sigma_0^\alpha + \sigma_1^\alpha}, \delta_0 + \delta_1 \right)$.
\end{enumerate}
\end{proposition}

Proposition \ref{prop:sum_CL_stable} provides a simple expression of $S_\alpha(\sigma,\beta,\delta)$ as a rescaling and shifting of an $S_\alpha(1,\beta,0)$. In fact given two independent random variables  $X\sim S_\alpha (\sigma,\beta,\delta)$ where $\sigma>0$ and $Z \sim S_\alpha (1,\beta,0)$, $X$ can  be expressed as
\begin{align}\label{eq:shift_rescal_simulation_stable}
X\stackrel{d}{=} \begin{cases} \sigma Z + \delta , &\text{ if } \alpha \neq 1 \\
\sigma Z + \delta + \frac{2}{\pi} \beta \sigma \log(\sigma) &\text{ if } \alpha=1.
\end{cases}
\end{align}

\emph{Tail behavior, moments and support:} 
The tail behavior of Gaussian distributions is well know; the survival function of $X\sim \mathcal{N}(0,1)$ asymptotically behaves like the following
$\Pro(X>x) \underset{x\rightarrow \infty}{\sim} \frac{1}{x\sqrt{2\pi}} e^{-x^2/2}.$
For general stable distributions ($\alpha \in (0,2)$), similar equivalents can be derived. The asymptotic behavior of a $X\sim S_\alpha(1,0,0)$ is of the order of a power function $x^{-\alpha}$,
which dominates the Gaussian tail equivalent, justifying the term 'heavy-tailed' distribution. 
\begin{theorem}[Theorem 1.2 \cite{nolan_univariate_2021}]\label{th:tail_stable}
Let $X\sim S_\alpha(\sigma,\beta,\delta)$ with $\alpha \in (0,2), \sigma>0, \delta \in \R$. Let $f_X$ be the corresponding density.
\begin{enumerate}
\item For $\beta \in (-1,1]$, $
\Pro(X>x) \underset{x\rightarrow \infty}{\sim} c_\alpha \sigma^\alpha(1+\beta)x^{-\alpha},$ and $
f_X(x) \underset{x\rightarrow \infty}{\sim} \alpha c_\alpha \sigma^\alpha (1+\beta) x^{-\alpha-1}$.

\item For $\beta \in [-1,1)$, $
\Pro(X<x) \underset{x\rightarrow -\infty}{\sim} c_\alpha \sigma^\alpha(1-\beta)x^{-\alpha},$ and $
f_X(x) \underset{x\rightarrow -\infty}{\sim} \alpha c_\alpha \sigma^\alpha (1-\beta) x^{-\alpha-1},$
\end{enumerate}
where $c_\alpha= \frac{1-\alpha}{2\Gamma(2-\alpha) \cos(\pi\alpha/2)}\ind_{\alpha\neq 1} + 
\pi^{-1} \ind_{ \alpha = 1}$.
\end{theorem}

Moments of any random variable can be expressed using the distribution function as $\E(|X|^p) = \int_{0}^\infty \Pro(|X|^p>x) dx$. For stable random variables, the existence of its moments is linked to the stability index $\alpha$ as a consequence of Theorem \ref{th:tail_stable}. Among $\alpha$-stable distributions ($\alpha \in (0,2]$), only the Gaussian distribution ($\alpha=2$) have finite variance (and moments of any order). Moreover only distributions with $\alpha \in (1,2]$ have a finite mean. More generally (see Proposition 1.2.16 in \cite{samoradnitsky_stable_2017})
\begin{align}\label{eq:moment_stable}
\E(|X|^p) &< \infty \quad \text{for all } p \in (0,\alpha),\\
\E(|X|^p) &= \infty \quad \text{for all } p \in [\alpha,+\infty).
\end{align}

Stable random variables either take values on the whole line $\R$ (in most cases) or in half lines when they are totally skewed ($\beta \in \{-1,1\}$). Let $X\sim S_\alpha(\sigma,\beta,\delta)$, the support of the associated density function $f_X$ is given by
\begin{align}\label{eq:support_stable}
supp(f_X) = \begin{cases}
[\delta- \sigma\tan(\pi\alpha/2), +\infty) & \alpha<1 \text{ and } \beta=1\\
(-\infty, \delta + \sigma\tan(\pi\alpha/2)]& \alpha<1 \text{ and } \beta=-1\\
\R & \text{otherwise}.
\end{cases}
\end{align}
A proof relying on an integral expression of the distribution function can be found in Section 3.2 of \cite{nolan_univariate_2021}. A parallel instructive approach involves a compound Poisson approximation and showing that $S_\alpha(\sigma,1,0)$ is the limit distribution of positive compound Poisson distributions (Proposition 1.2.11 of \cite{samoradnitsky_stable_2017}). 

\subsubsection*{Numerical evaluation of the stable densities}\label{sec:numerical_stable_densities}

A major holdback of the stable distribution lies in the lack of tractable formulas for the densities. In most cases, there is no explicit formulation of both the density and the distribution function of $S_\alpha(\sigma,\beta,\delta)$. However, some properties of these densities are known to provide numerical approximation algorithms. First, due to the exponential decay of their characteristic function, stable random variables exhibit . The law $\Pro_{X}$ is absolutely continuous, and its density is smooth. Indeed from \eqref{eq:def_stable_cf} any stable random variable $X\sim S_\alpha(\sigma,\beta,\delta)$ satisfies for $u\in \R$ and $k\in \N, |u^k\phi_X(u)| = |u|^k|\E\left(e^{iuX}\right)| = |u|^k e^{-\sigma^\alpha |u|^\alpha},$ providing a bound by an integrable function. Consequently, using Fourier theory, the density $f_X \in C^\infty(\R) \cap L^2(\R)$.

A first approach for numerical approximation of $f_X$ is to use Fourier inversion and numerical integration. In practice, directly inverting the expression \eqref{eq:def_stable_cf} of the characteristic function is numerically costly and can propagate errors due to the oscillating term and the infinite integration bound. For example, if $X\sim S_\alpha(1,\beta,0)$ the density of $X$ is given by
\begin{align*}
2\pi f_X(x) &= \int_{\R} e^{-iux}\phi_X(u)dx
= \begin{cases} 2\int_{0}^\infty e^{-u^\alpha} \cos(xu - \beta \tan(\frac{\pi\alpha}{2}) u^\alpha ) du \quad \text{if } \alpha \neq 1\\
2\int_{0}^\infty e^{-u} \cos\left(xu + \beta \frac2\pi u \log(|u|) \right) du \quad \text{if } \alpha=1.
\end{cases}
\end{align*}

Using a reformulation by Zolotarev of the stable densities and distribution functions (see \cite{zolotarev_one-dimensional_1986} and \cite{nolan_univariate_2021} for more details), more refined expressions can be derived, avoiding the numerical issues arising from the oscillatory nature of the integrand. First, using the reduction \eqref{eq:shift_rescal_simulation_stable}, one can focus on $Z\sim S_\alpha(1,\beta,0)$ 
and achieve the density of $X\sim S_\alpha(\sigma,\beta,\delta)$ by rescaling and shifting; $f_{X}(x) =\sigma^{-1}f_{Z}((x-\tilde\delta)/\sigma)$
where $\tilde\delta = \delta + (\frac{2}{\pi} \beta \sigma \log(\sigma))\ind_{\alpha=1}.$
In \cite{nolan_univariate_2021}, the following continuous function is introduced $V_{\alpha,\beta}$ in $(-\theta_0,\pi/2)$ where  $\theta_0 = \alpha^{-1} \arctan ( \beta \tan(\pi\alpha/2))\ind_{\alpha \neq 1} +
\frac\pi2 \ind_{\alpha=1}$,
\begin{align*}
V_{\alpha,\beta}(\theta) = \begin{cases}
\left(\cos(\alpha \theta_0\right)^{\frac{1}{\alpha-1}} \left( \frac{\cos(\theta)}{\sin(\alpha(\theta_0 + \theta))} \right)^{\frac{\alpha}{\alpha-1}} \frac{\cos(\alpha \theta_0 + (\alpha-1)\theta}{\cos(\theta)} & \alpha \neq 1\\
\frac{2}{\pi} \left( \frac{\frac\pi2 + \beta \theta}{\cos(\theta)} \right) \exp \left(\beta^{-1} (\frac\pi2 + \beta \theta) \tan(\theta) \right) & \alpha=1, \beta \neq 0.
\end{cases}
\end{align*}
Note that for $\alpha=1$ and $\beta=0$, the function $V_{\alpha,\beta}$ is not defined but in this case, $Z$ is a Cauchy random variable for which a simple formula of its density is available. The most important feature of $V$ is its non oscillatory behavior over the integration interval $(-\theta_0,\pi/2)$ (Lemma 3.9 in \cite{nolan_univariate_2021}). The following formulas for the density are derived.

\begin{theorem}[Theorem 3.3 of \cite{Nolan1997}]\label{th:density_stable_expression}
Let $Z \sim S_\alpha(1,\beta,0)$ and denote by $f_{Z}$  the corresponding density function. Then, for $\alpha \neq 1$ and  $\theta_0 = \alpha^{-1} \arctan ( \beta \tan(\pi\alpha/2))\ind_{\alpha \neq 1} +
\frac\pi2 \ind_{\alpha=1}$, the density $f_Z$ can be expressed for $x\in \R$ as 
\begin{align*}
f_{Z}(x)&=\begin{cases} \frac{\alpha x^{\frac{1}{\alpha-1}}}{\pi |\alpha-1|} \displaystyle \int_{\theta_0}^{\pi/2} V_{\alpha,\beta}(\theta) \exp\left(-x^{\frac{\alpha}{\alpha-1}} V_{\alpha,\beta}(\theta) \right) d\theta & x>0,\\
\pi^{-1}\Gamma(1+ \alpha^{-1}) \cos(\theta_0) (\cos(\alpha\theta_0))^{1/\alpha} & x=0,\\
f_{-Z,\alpha}(x) & x<0.
\end{cases}
\end{align*}
where $-Z \sim S_\alpha(1,-\beta,0)$. For $\alpha=1$, the density $f_Z$ can be expressed for $x\in \R$ as 
\begin{align*}
f_{Z}(x)&=\begin{cases} \frac{1}{2\beta} e^{-\frac{\pi x}{2\beta}}\displaystyle \int_{-\pi/2}^{\pi/2} V_{1,\beta}(\theta) \exp\left(-e^{-\frac{\pi x}{2\beta}} V_{1,\beta}(\theta) \right) d\theta & \beta \neq 0,\\
\frac{1}{\pi(1+x^2)} & \beta=0.
\end{cases}
\end{align*}
\end{theorem}
Combining \eqref{eq:shift_rescal_simulation_stable} and Theorem \ref{th:density_stable_expression} lead to more robust numerical computations of the aforementioned densities. More details on the numerical approximation of stable densities and distribution functions can be found in \cite{Nolan1997}.
The following figures display densities of $S_\alpha(\sigma,\beta,\delta)$ with varying parameters $(\alpha,\sigma,\beta,\delta)$ and illustrate their effect on the stable densities.

\subsection{The $\alpha$-stable Lévy process}\label{sec:stable_levy_process}
%\textbf{References:}\cite{samoradnitsky_stable_2017,sato_levy_1999,cont_financial_2003}
The infinite divisibility of stable random variables provides via the Lévy-Khintchine theorem a useful expression of the characteristic function  (see Remark \ref{rem:infinite_divisibility}). It also establishes a connection between stable random variables and Lévy processes. Specifically, for every Lévy process $(X_t)_{t\in \R_+}$ the increment $X_t$ is infinitely divisible. Conversely, if $Y$ is an infinitely divisible random variable, it is always possible to construct a Lévy process $(X_t)_{t\in \R+}$ such that $X_1 \stackrel{d}{=} Y$ (Theorem 7.10 of \cite{sato_levy_1999}). An $\alpha$-stable Lévy process is then defined as a Lévy process $(X_t)_{t\geq 0}$ such that $X_1$ is a stable random variable with stability index $\alpha \in (0,2]$. Every Lévy process $(X_t)_{t\in \R_+}$ is characterized by the Lévy triplet of the infinitely divisible random variable at time $1$.
\begin{remark}
For $\alpha \in (0,2)$, the Lévy triplet of $(X_t)_{t\in \R_+}$ is $(b,0,\nu)$ where $b\in \R$ is the drift term, the volatility parameter is null, and the Lévy measure is given by
 \begin{align}\label{eq:def_levy_measure_stableProcess}
\nu(dx) = \frac{P}{x^{1+\alpha}} \ind_{x>0} + \frac{Q}{|x|^{1+\alpha}}\ind_{x<0},
\end{align}
where $P,Q\geq 0$ such that $P+Q>0$. In this case, the total mass is infinite, i.e. $\nu(\R)= \infty$. The $\alpha$-stable Lévy process is an important example of a Lévy process with infinite jump activity.
For $\alpha=2$, the process $(X_t)_{t\in \R_+}$ is a rescaled Brownian motion of Lévy triplet $(b,\sigma,0)$.
\end{remark}

\subsubsection*{Selfsimilarity and stable parameters}

The stability of $X_1$ transfers to all marginals $X_t$ for $t>0$, thanks again to the infinite divisibility property. Let $X_1 \sim S_\alpha(\sigma,\beta,\delta)$. For all $t>0$ the characteristic function can be expressed as $\phi_{X_t}(u)= e^{t\eta(u)}$ where $\eta \in \R \rightarrow \C$ is usually called the Lévy symbol (see \cite{sato_levy_1999}). Consequently, it appears that for $u\in \R, \phi_{X_t}(u) = \phi_{X_1}(u)^t.$
From the general expression of the stable characteristic function \eqref{eq:def_stable_cf} $\forall t>0$, $X_t$ is an $\alpha$- stable random variable
\begin{align*}
\phi_{X_t}(u) = \phi_{X_1}(u)^t = \begin{cases}
    \exp \left( - \left( t^{1/\alpha}\sigma\right)^\alpha |u|^\alpha \left(1- i\beta sign(u) \tan \left(\frac{\pi\alpha}{2} \right) \right) + i\delta t u \right) &\text{ if } \alpha \neq 1, \\
    \exp \left(- t\sigma |u|(1+ i\beta \frac{2}{\pi} sign(u) \log(|u|)) + i \delta t u \right) & \text{ if } \alpha =1.
\end{cases}
\end{align*}
Moreover, as a stable random variable, the parameters of $X_t$ are merely computed by rescaling those of $X_1$; i.e. for all $t>0$, $X_t \sim S_\alpha \left(t^{1/\alpha} \sigma,\beta, t\delta\right).$ It leads to the following equality in law for $X_t$
\begin{align}\label{eq:X_t_parameter_stable}
    X_t \stackrel{d}{=}  \begin{cases}
        t^{1/\alpha}X_1 + \left(t-t^{1/\alpha} \right)\delta, &\text{ if } \alpha \neq 1, \\
        t X_1 + \frac{2}{\pi}\beta t \log(t),& \text{ if } \alpha=1.
    \end{cases}
\end{align}  
This property is known as broad-sense selfsimilarity. Selfsimilar processes are an example of processes that exhibit an invariance property connecting a scaling in time to a scaling in space (see \cite{sato_levy_1999,samoradnitsky_stable_2017} for more details). A stochastic process $(X_t)_{t\geq 0}$ is said to be broad-sense selfsimilar with Hurst index $H>0$ (Definition 13.4 \cite{sato_levy_1999}) if, for any $r>0$ there exist a function $b: (0,\infty) \rightarrow \R$ such that $$ \{ X_{rt}, t\geq 0 \} \stackrel{d}{=} \{ r^{H} X_t + b(t), t\geq 0 \}.$$ For example, any drifted Brownian motion $(\sigma W_t + bt)_{t\in \R_+}$ with $\sigma>0$ and $b\in \R$ is a broad-sense selfsimilar process with Hurst index $1/2$. From the fact that $\phi_{X_t}(u) = \phi_{X_1}(u)^t$ it is straightforward to see that the only broad-sense selfsimilar Lévy process are the $\alpha$-stable processes. Morally, the $\alpha$-stable Lévy process can be seen as an extension of the Brownian motion through its selfsimilarity property.

\begin{remark}(From the Lévy parameters $(P,Q,\alpha)$ to the stability parameters $(\alpha,\sigma,\beta,\delta)$).
Given an $\alpha$-stable Lévy process, our access to its characteristic parameters may be limited to the Lévy triplet. For this reason, it is important to explicitly provide formulas connecting the ($\alpha,\sigma,\beta,\delta$) parameterization of stable distributions to the Lévy parameters ($\alpha,P,Q$). Let $X$ be an $\alpha$-stable Lévy process with $\alpha \in (0,2)$ and with Lévy triplet $(0,0,\nu)$ where $\nu$ is given by \eqref{eq:def_levy_measure_stableProcess}. Let $X_1 \sim S_\alpha(\sigma,\beta, \delta)$ for $\sigma>0, \beta \in [-1,1]$ and $\delta \in \R$ then the following formulas can be derived.
\begin{align}\label{eq:parameter_levy-to-stable}
\sigma= \begin{cases}
        \left(  \frac{(P+Q)}{\alpha} \Gamma(1-\alpha) \cos\left(\frac{\pi\alpha}{2}\right))\right)^{1/\alpha} &  \alpha\neq 1\\
        \frac{\pi}{2}(P+Q),&  \alpha=1
        \end{cases}, & 
        & \beta=\frac{P-Q}{P+Q}, &
        &\delta= \begin{cases} \frac{Q-P}{1-\alpha} & \alpha \neq 1 \\  c(P-Q)&  \alpha=1
        \end{cases}
\end{align}
where $c=\int_{1}^\infty  \frac{sin(r)}{r^2}dr + \int_{0}^1 \frac{sin(r)-r}{r^2}dr.$ It comes from a straightforward computation using the Lévy-Khintchine formula \eqref{eq:levy_khintchine} with $\sigma=\delta=0$ and the complex integral from Lemma 14.11 of \cite{sato_levy_1999}. The drift term $\delta$ depends on the convention form of the Lévy-Kintchine formula (see \cite{sato_levy_1999}). 
\end{remark}

\subsection{Trajectory simulation}
This section is organized as follows. A sampling algorithm for a general $X\sim S_\alpha(\sigma,\beta,\delta)$ random variable is presented: Algorithm  \ref{algo:simulation_S(sigma,beta,delta)}. It is based on a Box-Muller-type algorithm proposed by Chamber and Mallow in \cite{chambers1976method}. Then for an $\alpha$-stable Lévy process, using the independence and stationarity of its increments, a sample of the trajectory can be computed via Algorithm \ref{algo:simulation_Stable_skeleton)}. Several references are devoted to simulating stable random variables and processes \cite{cont_financial_2003,weron2005computer}.
To assess the accuracy of these algorithms, the increments histogram is compared to the stable densities which are numerically computed based on Theorem \ref{th:density_stable_expression}.
\subsubsection{Sampling from stable distributions}\label{sec:simulation_stable}

The $2-$stable distribution is a Gaussian $X\sim \mathcal{N}(b,A^2)$, and the Box-Muller algorithm is widely used for simulating $X$. Let $U_1,U_2$ be two independent uniform random variables on $[0,1]$, then $$b + A\sqrt{-2\log(U_1)} \cos(2\pi U_2) \sim \mathcal{N}(b,A^2).$$ An analogous algorithm can be devised for simulating the Cauchy distribution with parameter $X\sim S_1(\sigma,0,\delta)$ as $X \stackrel{d}{=} \sigma \tan \left( U \right) +\delta$,
where $U\sim \mathcal{U}\left(-\frac\pi2,\frac\pi2 \right)$.

In the more general setting of $S_\alpha(\sigma,\beta,\delta)$ simulation simplifies through rescaling and shifting to that of $S_\alpha (1,\beta,0)$ \eqref{eq:shift_rescal_simulation_stable}. The following Theorem gives a general method to sample from $S_\alpha(1,\beta,0)$, for all $\alpha \in (0,2)$ . 
\begin{theorem}[Theorem 1.3 \cite{nolan_univariate_2021}]\label{Th:simulation_stable}
    Let $\alpha \in (0,2)$ and  $U,V$ be two independent random variables such that $U\sim \mathcal{U}\left([-\frac\pi2, \frac\pi2] \right)$ and $V\sim \mathcal{E}(1).$
Let $\theta = \alpha^{-1}\arctan\left(\beta \tan\left(\frac{\pi\alpha}{2} \right) \right)$ for  $\beta \in [-1,1]$  and $\alpha \neq 1$ then
 \begin{align}\label{eq:box_muller_stable}
     X \sim \begin{cases}
\frac{\sin(\alpha(\theta + U))}{\left(\cos(\alpha \theta)\cos(U)\right)^{1/\alpha}} \left(\frac{\cos(\alpha \theta + (\alpha-1)U)}{V}\right)^{(1-\alpha)/\alpha} & \text{ if } \alpha \neq 1, \\
  \frac{2}{\pi} \left[\left( \frac{\pi}{2} + \beta U\right)\tan(U) - \beta \log\left( \frac{\frac{\pi}{2}V\cos(U)}{\frac{\pi}{2}+ \beta U} \right) \right]       & \text { if } \alpha=1,
     \end{cases}
 \end{align}
 is distributed as $S_\alpha (1,\beta, 0).$
\end{theorem}
\begin{remark}
In particular, for symmetric stable distribution (i.e $\beta=0$), the formula simplifies in
\begin{align*}
    \frac{\sin(\alpha U)}{\cos(U)^{1/\alpha}} \left[\frac{\cos((\alpha-1)U)}{V} \right]^{\frac{1-\alpha}{\alpha}} &\sim S_\alpha(1,0,0),\\
    \tan(U) &\sim S_1(1,0,0),
\end{align*}
which aligns with the standard Cauchy simulation algorithm.
\end{remark}
The following Proposition establishes that any stable random variable $S_\alpha(1,\beta,0)$ can be expressed as a linear combination of positively skewed stable random variables $S_\alpha(1,1,0)$. These simple distributions relate to the class of stable subordinator processes that are key to sample of CTS Lévy processes (Section \ref{sec:simulation_tempered_stable}).
\begin{proposition}[Proposition 1.2.13 of \cite{samoradnitsky_stable_2017}]\label{prop:stable=CL_skewed}
Let $X\sim S_\alpha (1,\beta,0)$ where $\alpha \in (0,2)$ and $\beta \in [-1,1]$. For $Y_1,Y_2$ are independent random variables distributed as $S_\alpha (1,1,0)$, it holds that
\begin{align}\label{eq:skewed_simu}
X \stackrel{d}{=} \begin{cases}
    \left( \frac{1+\beta}{2} \right)^{1/\alpha} Y_1 - \left( \frac{1-\beta}{2} \right)^{1/\alpha} Y_2& \text{ if } \alpha\neq 1,\\
     \left( \frac{1+\beta}{2} \right)Y_1 - \left( \frac{1-\beta}{2} \right) Y_2 + \sigma\left( \frac{1+\beta}{\pi} \log\left( \frac{1+\beta}{2}\right)-\frac{1-\beta}{\pi} \log\left( \frac{1-\beta}{2}\right)\right)& \text{ if } \alpha=1.
\end{cases}
\end{align}
\end{proposition}

\begin{remark}
An arbitrary stable distribution $S_\alpha(\sigma,\beta,\delta)$ can be sampled through the combined application of Theorem \ref{Th:simulation_stable} and \eqref{eq:shift_rescal_simulation_stable}. An alternative approach involves using Theorem \ref{Th:simulation_stable} to simulate $2$ copies of an $S_\alpha (1,1,0)$ and subsequently, by Proposition \ref{prop:stable=CL_skewed}, obtaining realizations from $S_\alpha(1,\beta,0)$ which can be shifted and rescaled  (as in \eqref{eq:shift_rescal_simulation_stable}) into a realization of $S_\alpha(\alpha,\beta,\delta)$. Algorithm \ref{algo:simulation_S(sigma,beta,delta)} is derived for the former approach. 
\end{remark}

\subsubsection{Simulating trajectories of $\alpha$-stable Lévy processes}
Among the rich class of Lévy processes, the $\alpha$-stable Lévy process is an example where exact simulation of the increments is achievable (see \cite{cont_financial_2003}). In practice, simulation methods typically focus on constructing a discrete skeleton of the Lévy process over a fixed grid. Let $X$ be a Lévy process, and consider a fixed time horizon $T>0$ and $n\geq 1$ the number of observations. The observation or sampling rate is defined as  $\Delta= \frac{T}{n}$.
\begin{itemize}
\item Leveraging the fact that the increments are independent and distributed as the first one $X_\Delta$, trajectory simulation reduces to sampling $n$ independent copies of $X_\Delta$.

\item Using the (broad-sense) selfsimilarity property (see Section \eqref{sec:stable_levy_process}) it is sufficient to draw $n$ independent copies $(Y_i)_{i=1}^n$ of $X_1$ and use \eqref{eq:X_t_parameter_stable}. To sample $X_1$ given the Lévy parameters $(\alpha,P, Q)$ a translation into the stable parameters $(\alpha,\sigma,\beta,\delta)$ is achieved using the formula devised in \eqref{eq:parameter_levy-to-stable}.

\item  Then by summing the increments we derive a discrete time numerical approximation of the continuous time Lévy process trajectory. Usually for displaying, the value of $X_t$ at time $t \in [0,T]$ where $t\neq k\Delta$ for $k\leq n$ is computed via a linear interpolation. For $t\in [0,T]$ we define

\begin{align}\label{eq:skeleton_approximation}
\widetilde{X}_t = \begin{cases}
\sum_{j=0}^{k} \Delta^{(j)} X & \text{if } t=k\Delta, k = 0,\cdots, n\\
\text{ linearly interpolated, }  & \text{elsewhere.}
\end{cases}
\end{align}
\end{itemize}

\begin{remark}\label{rem:self-similarity_computational_cost}
In terms of complexity, the sampling algorithms are equivalently costly for low and high-frequency observations. This is mostly due to the selfsimilarity property which reduces the sampling of $X_\Delta$ to that of $X_1$.
\end{remark}
\begin{center}
\begin{algorithm}
\textbf{Step 1. }Generate independently $U\sim \mathcal{U}[-1,1]$ and $V\sim \mathcal{E}(1)$.

\textbf{Step 2: }\If{$\alpha=1$}
{
$Z=\frac{2}{\pi} \left[\left( \frac{\pi}{2} + \beta U\right)\tan(U) - \beta \log\left( \frac{\frac{\pi}{2}V\cos(U)}{\frac{\pi}{2}+ \beta U} \right) \right] $. \\
\KwRet{ $\sigma Z + \delta+ \frac{2}{\pi}\beta \sigma \log(\sigma)$}
}
\Else
{
$Z= \frac{\sin(\alpha(\theta + U)}{\left(\cos(\alpha \theta)\cos(U)\right)^{1/\alpha}} \left(\frac{\cos(\alpha \theta + (\alpha-1)U)}{V}\right)^{(1-\alpha)/\alpha} $. \\
\KwRet{ $\sigma Z + \delta$.}
}
\caption{Simulation of $X\sim S_\alpha(\sigma,\beta,\delta)$}
\label{algo:simulation_S(sigma,beta,delta)}
\end{algorithm}
\begin{algorithm}[H]
\textbf{Step 1. } Compute $(\sigma,\beta,\delta)$ from $(\alpha,P,Q,b)$ using \eqref{eq:parameter_levy-to-stable}.

\textbf{Step 2. } Sample $n$ independent copies $(Y_{i,1})_{i=1}^n$ of $X_1$ using Algorithm \ref{algo:simulation_S(sigma,beta,delta)} where $Y_{0,1}=0$.

\textbf{Step 3. } Compute $Y_{i,\Delta} =  \Delta^{1/\alpha} Y_{i,1} + (\Delta-\Delta^{1/\alpha}) \delta + \frac2\pi\beta \Delta \log(\Delta)\ind_{\alpha=1}$ for $i \in \{1,\cdots, n\}$.

\KwRet{$\left(\sum_{i=1}^\ell Y_{i,\Delta}\right)_{\ell=0}^{n}$.}

\caption{Sampling $(\widetilde X_t)_{t\in \R_+}$ on $\{0,\Delta,\cdots,n\Delta\}$ for an $\alpha$-stable Lévy process with triplet $(b,0,\nu)$ ($\nu$ given by \eqref{eq:stable_levymeasure})}
\label{algo:simulation_Stable_skeleton)}
\end{algorithm}
\end{center}

\begin{figure}
\begin{center}
\includegraphics[scale=0.34]{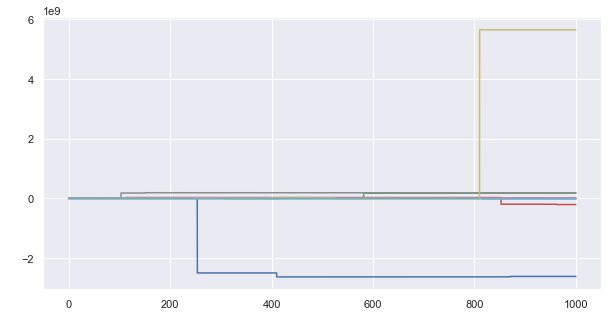}
\includegraphics[scale=0.44]{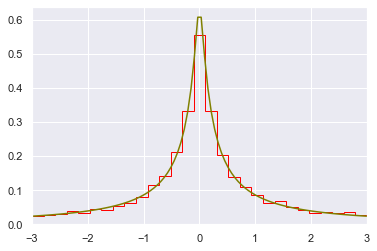}\\
\includegraphics[scale=0.34]
{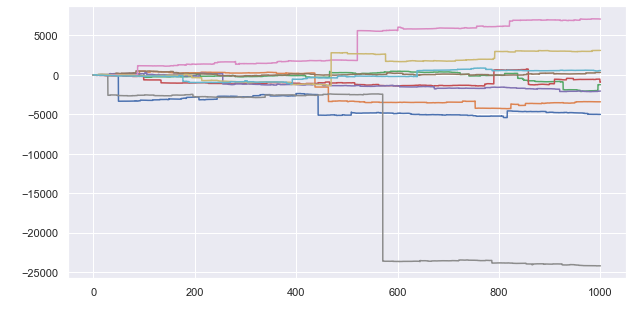}
\includegraphics[scale=0.44]{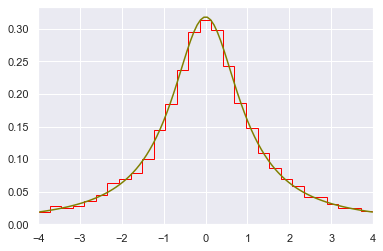}\\

\includegraphics[scale=0.34]{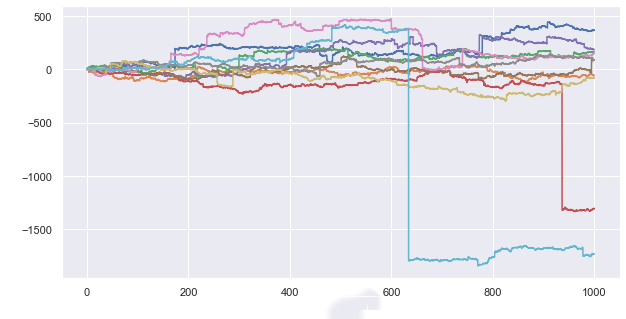}
\includegraphics[scale=0.42]{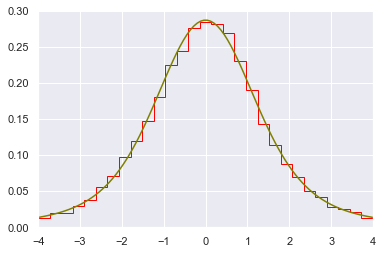}
\caption{Trajectories of a symmetric $\alpha$- stable Lévy process with Lévy triplet $(0,0,|x|^{-1-\alpha})$ where $\alpha \in \{0.5,1,1.5\}$, $\Delta=1$, $n=1000$. Histogram of $X_\Delta$ and the stable density $f_{X_\Delta}$ are displayed. }
\label{fig:trajectories_stable_process}

\end{center}
\end{figure}

\section{Classical Tempered Stable Lévy processes (CTS)}

\subsection{Tempering an $\alpha$-stable process}
%\textbf{References:}\cite{blumenthal_sample_1961,rosinski_generalized_2010,carr_fine_2002,kuchler_tempered_2013}
Lévy processes can be constructed using various transformations of other Lévy processes e.g. by linear transformation or subordination (see \cite{cont_financial_2003}). Tempering is one such procedure. It relies on multiplying a Lévy measure $\nu$ of a Lévy process $(X_t)_{t\in \R_+}$ by a decreasing exponential function on each half-axis. This results in a modified measure $$\tilde \nu (dx) = (e^{-Ax}Px^{-1-\alpha} \ind_{x>0} + e^{-B|x|}Q|x|^{-1-\alpha}\ind_{x<0}),$$ which still satisfies the properties of a Lévy measure. The associated Lévy process $(\tilde X_t)_{t\in \R_+}$ governed by $\tilde\nu$ is referred to as the Esscher transform \cite{Esscher1932} (see Chapter 9 of \cite{cont_financial_2003}). This transformation has proven useful in actuarial science \cite{gerber1993option}. Furthermore, exponential tilting is also strongly connected with Monte-Carlo methods, including rejection sampling and importance sampling (\cite{Baeumer2010TemperedSL,asmussen2007stochastic}). These connections play a critical role in the sampling algorithm for the tempered processes (see Section \ref{sec:simulation_tempered_stable}).

Specifically, tempering an $\alpha$-stable Lévy process leads to the Classical Tempered Stable process (CTS), defined as a Lévy process $(X_t)_{t\in \R_+}$ with a Lévy triplet $(b,0,\nu)$ where $b\in \R$ and \begin{align}\label{eq:levy_measure_tempered_stable}
\frac{\nu(dx)}{dx} = \frac{Pe^{-Ax}}{x^{1+\alpha}} \ind_{x>0} + \frac{Qe^{-B|x|}}{|x|^{1+\alpha}} \ind_{x<0}.
\end{align}
where $P,Q,A,B \in \R_+$ such that $P+Q>0$ and $A+B>0$, $\alpha \in (0,2)$. When $P,Q \neq 0$ the process is referred to as \emph{bilateral CTS} (see \cite{kawai_simulation_2011}) encompassing the totally positively (resp. negatively) skewed CTS process $Q=0$ (resp. $P=0$). 
\begin{remark}
Tempering can be interpreted through its effect on the trajectories $t\rightarrow X_t$. The exponential factor diminishes the intensity of the big jumps while preserving the stable nature of small jumps. As a result of tempering the Lévy measure, the activity index $\alpha$ can take negative values as $\int_{\R}(x^2\wedge 1)\nu(dx) <\infty$ for $\alpha \in (-\infty,2)$. Indeed, if $\alpha<0$, the total Lévy measure is finite $\nu(\R) = (PA^\alpha + QB^\alpha) \Gamma(-\alpha)<\infty$. In this case, the process $(X_t)_{t\in \R_+}$ is a Compound Poisson process offering finite jump activity models with intricate structures.
 If $\alpha \in [0,1)$, the process has infinite activity since $\nu(\R)= \infty$ and a finite variation $\int_{|x|<1} |x| \nu(dx) <\infty$. Conversely, if $\alpha \in [1,2)$, the total variation is infinite, i.e. $\int_{|x|\leq 1} |x|\nu(dx) = \infty$. 
\end{remark}

\subsubsection{Distributional properties with a view toward simulation}

In this Section, we select the key tools from the many properties of CTS distributions, to tackle the sampling problem. For further details, we refer to \cite{kuchler_Tappe_2008,cont_financial_2003,rosinski_generalized_2010}.
First, we highlight that closed parametric formulas for their characteristic functions \eqref{eq:cf_tempered} are available. These computations rely on the Lévy-Khintchine formula, but where the classical small and big jump decomposition is dropped. In particular, let $(X_t)_{t\in \R_+}$ to be a CTS Lévy process with triplet $(b_\nu, 0,\nu)$ where $b_\nu$ is a drift term $b_\nu =\int_{|x|>1}x \nu(dx) = PA^{\alpha-1}\Gamma(1-\alpha,A) - QB^{\alpha-1}\Gamma(1-\alpha,B)$ ensures that for $t>0$
\begin{align}\label{eq:levy_khintchine_tempered_general}
\log(\phi_{X_t}(u))&= itub_\nu + t\int_{\R}(e^{iux}-1 - iux \ind_{[-1,1]}) \nu(dx)= t\int_{\R}(e^{iux} - 1 - iux ) \nu(dx).
\end{align}
Using complex contour integration \cite{kuchler_Tappe_2008} or with a power series argument \cite{cont_financial_2003} the following result can be derived
\begin{proposition}[Proposition 4.2 of \cite{cont_financial_2003}]
Let $(X_t)_{t\in \R}$ be a CTS Lévy process with Lévy triplet $\left(b_\nu, 0,\nu\right)$ where $\nu$ defined as in \eqref{eq:levy_measure_tempered_stable} and $b_\nu=\int_{|x|>1}x \nu(dx)$. For $t\in \R_+$, the characteristic function exponent can be expressed for $u\in \R$ as 
\begin{align}\label{eq:cf_tempered}
\log(\phi_{X_t}(u))= \begin{cases} 
tP\Gamma(-\alpha)\left((A-iu)^\alpha - A^\alpha + iu\alpha A^{\alpha-1} \right) + tQ\Gamma(-\alpha)\left((B+iu)^\alpha - B^\alpha + iu\alpha B^{\alpha-1}  \right) & \alpha \neq 1,\\
iut(P-Q) + tP(A-iu)\log\left(1-\frac{iu}{A}\right) + tQ(B-iu)\log\left(1-\frac{iu}{B}\right)& \alpha = 1.
\end{cases}
\end{align}
\end{proposition}
For $\alpha>1$ and $A=B=0$ the exponent \eqref{eq:cf_tempered} yields an equivalent formulation of the $\alpha$-stable characteristic function exponent \eqref{eq:def_stable_cf}. For instance for $\alpha \in (1,2)$ and $A=B=0$ and $Q=0$ the characteristic function exponent is given by $$\log(\phi_{X_t}(u)) = tP\Gamma(-\alpha)(-iu)^\alpha = -tP\alpha^{-1}\Gamma(1-\alpha)\cos(\pi\alpha/2) |u|^\alpha (1- i sign(u)\tan(\pi\alpha/2)),$$
aligning with Formula \eqref{eq:parameter_levy-to-stable}.
\begin{remark}
For $\alpha$-stable Lévy processes the self-similarity property reduces the computational cost (see Remark \ref{rem:self-similarity_computational_cost}). However, the characteristic function expression \eqref{eq:cf_tempered} shows that this property does not extend to the CTS processes. Sampling increments in low-frequency (small $\Delta$) is computationally more expensive than for high-frequency (large $\Delta$) (see Section \ref{sec:simulation_tempered_stable}).

One of the key property of CTS processes is that, unlike $\alpha$-stable process, they have finite moments. In fact, $\E(|X_t|^\beta)<\infty$ for all $\beta>0$ (Proposition 2.7 \cite{rosinski_generalized_2010}). Furthermore, they also have finite exponential moments: $\E(e^{\theta X})<\infty$ if $\theta \leq \min(A,B)$ ((iv) of Proposition \cite{rosinski_generalized_2010}). In the general setting, the cumulant generating function can be derived (see Remark 2.8 of \cite{kuchler_tempered_2013}). It provides parametric estimators of the Lévy measure parameters $\alpha,P,A,Q$ and $B$ (see Section 6 of \cite{kuchler_Tappe_2008} for more details). 
\end{remark}

\begin{remark}(\textit{Long and short time behaviors})
The asymptotic behavior of CTS processes resembles a stable Lévy process in a small time scale while in a long time scale, it approximates a Brownian motion. After appropriate shifting and rescaling, the following convergence in distribution holds
\begin{align*}
    \left(\frac{X_{\Delta t}-b_t(\Delta)}{\Delta^{1/\alpha}} \right)_{t\in \R_+} \underset{\Delta \rightarrow 0}{\rightarrow} (S_t)_{t\in \R_+} & & \left(\frac{X_\Delta t - b_t(\Delta)}{\Delta^{1/2}}\right)_{t \in \R_+} \underset{\Delta \rightarrow 0}{\rightarrow} (W_t)_{t\in \R_+},
\end{align*}
where $S$ is a strictly stable $\alpha$-stable process and $W$ is a drifted Brownian motion, $b_t(h)$ is a adequate shifting term (see Theorem 3.1 of \cite{rosinski_generalized_2010} for more details). 

\end{remark}
\subsubsection{CTS densities and numerical approximation}\label{sec:tempered_stable_density}
The class of CTS process stands among the Lévy process for which the smoothness of the density is not a time-dependent property (see \cite{sato_levy_1999}). Indeed, it satisfies the Orey's smoothness criterion for infinitely divisible distributions $\liminf_{r\rightarrow 0} \frac{\int_{|x|\leq r}{x^2 \nu(dx)}}{x^{2-\alpha}}>0$ (see Proposition 28.3 \cite{sato_levy_1999} and \cite{Orey68}). This means that for all $t>0$, $X_t$ has a density $f_{X_t}$ which is of class $C^\infty$. 
In the finite variation case $\alpha \in (0,1)$, an extensive study of the unimodality and asymptotic behavior of $f_{X_t}$ is provided in \cite{kuchler_tempered_2013}. In particular, Theorem 7.10 of \cite{kuchler_tempered_2013} states that for $t>0$ there exists a constant $C>0$ such that, $$f_{X_t}(x) \underset{x\rightarrow +\infty}{\sim} C \frac{e^{-Ax}}{x^{1+\alpha}}.$$
Compared to the asymptotic equivalent for the $\alpha$-stable density (Theorem \ref{th:tail_stable}), the tempering adds a decreasing exponential factor. This ensures finite moments but still preserves a quasi-stable behavior in a short time.

For the numerical aspects of CTS densities $f_{X_t}$, a natural approach is to use the Fourier inverse formula leveraging the explicit characteristic function \eqref{eq:cf_tempered}. To the best of our knowledge, no simplified formulation as in Section \ref{sec:numerical_stable_densities} can be used to avoid the oscillating effect of the characteristic function. Morally one can expect those to be tempered by the exponential term leading to a more robust integration method. We denote by $\tilde{f}_{X_t,M}$ the following approximation function of $f_{X_t}$ defined for $x\in \R$ as $$\tilde{f}_{X_t,M}(x) = (2\pi)^{-1} \int_{-M}^M \phi_{X_t}(u) e^{-iux}du.$$
The choice of $M>0$ is crucial and can be numerically costly. In practice one might select $M$ large enough such that $\tilde{f}_{X_t,M}$ does not change, i.e. for $\eta$ small enough $|\tilde{f}_{X_t,M}(x)-\tilde{f}_{X_t,M'}(x)|\le \eta$ for all $x$ and any  $M'\geq M$.

\begin{remark}\label{rmk:totally_skewed_stable}
For totally positively skewed CTS processes $(Q=0)$, the density can be related to that of a totally skewed $\alpha$-stable process ($\beta=1$) for which efficient numerical algorithms are available (see Section \ref{sec:numerical_stable_densities}).  Let $(Y_t^+)$ be a centered and totally positively skewed CTS distribution i.e. whose characteristic function is given by 
$\E(e^{iuY_t^+})= \exp \left(t \int_{\R_+}(e^{iux} - 1 - iux) \frac{Pe^{-Ax}}{x^{1+\alpha}} dx \right).$ We denote by $(S_t^+)_{t\in \R_+}$ the totally positively skewed $\alpha$-stable Lévy process ($A=0$) such that $S_t^+ \sim S_\alpha(t\sigma,1,0)$ where $\sigma^\alpha = P\Gamma(1-\alpha)\alpha^{-1}\cos(\pi\alpha/2)\ind_{\alpha\neq 1} + \frac{P}{2\pi}\ind_{\alpha=1}$. 
From Proposition 1 of \cite{meerschaert_tempered_2014} the density of $Y_t^+$ can be expressed as the following transformation of that of $S_t^+$,
\begin{align}\label{eq:density_CTS_stable}
f_{Y_t^+}(x) = \begin{cases}
e^{-Ax -(1-\alpha) tP\Gamma(-\alpha)A^\alpha} f_{S_t^+}\left(x - tP\alpha\Gamma(-\alpha)A^{\alpha-1} \right) & \alpha \neq 1,\\
e^{-Ax +tPA} f_{S_t^+}(x - Pt (1+ \log(A)) ) &\alpha =1.
\end{cases}
\end{align}

\end{remark}

\subsection{Trajectory simulation algorithms}\label{sec:simulation_tempered_stable}

Among the many algorithms designed for sampling CTS distributions and processes on $\R$, the Bauemer-Merchaer acceptance-rejection algorithm \cite{Baeumer2010TemperedSL} strikes a good balance between simple implementation, parameter tuning, and computing efficiency. It relies on the simulation and rejection of stable increments for which efficient algorithms are established (see Section \ref{sec:simulation_stable}). A complete comparison between the Bauemer-Merchaer method and other algorithms (such as Devroye's algorithm \cite{devroye2009random}, Compound Poisson approximation \cite{cont_financial_2003}, shot noise series approximation \cite{rosinski_series_2001}) is discussed in \cite{kawai_simulation_2011} which concludes to its numerical efficiency in small time. 
 In this section, we present the complete bilateral case. For simplicity, we assume that the index $\alpha\in (0,1)\cup(1,2)$, as similar computations can be conducted when $\alpha=1$. First reduction to totally skewed CTS processes is discussed. Then, leveraging Remark  \ref{rmk:totally_skewed_stable}, the Bauemer-Merchaer algorithm is outlined (Algorithms \ref{algo:simulation_Y_t^+} and \ref{algo:simulation_Tempered_skeleton)}). A computational time analysis is conducted in both finite and infinite variation cases. Finally, trajectories are displayed (Figure \ref{fig:CTS_trajectories_histograms}), and a comparison between the histogram of the samples of $X_\Delta$ and the numerical computation of the CTS densities (as discussed in Section \ref{sec:tempered_stable_density}) evaluates the accuracy of the procedure.

\subsubsection*{Reduction to totally positively skewed CTS processes}

Similarly to $\alpha$-stable processes (see Proposition \ref{prop:stable=CL_skewed}, simulating general CTS process with bilateral Lévy measure $\nu(dx)$ as in \eqref{eq:levy_measure_tempered_stable} reduces to sampling two totally positively skewed CTS processes. In the following, $(X_t)_{t\in \R_+}$ is a CTS process with Lévy triplet $(b_\nu,0,\nu)$ and $b_\nu = \int_{|x|>1}x \nu(dx)$. The characteristic function exponent \eqref{eq:cf_tempered} can be separated into two components. Consequently, one can decompose $X_t = Y_t^{+} - Y_{t}^{-}$ where $Y^+$ and $Y^{-}$ are independent and their characteristic function is given for $t>0$ by
\begin{align}\label{eq:tempered_subordinator_cf}
\E(e^{iuY_t^+}) &=  e^{tP\Gamma(-\alpha) \left( (A-iu)^\alpha-A^\alpha + iu\alpha A^{\alpha-1} \right)}\\
\E(e^{iuY_t^-}) &= e^{{tQ\Gamma(-\alpha) \left( (B-iu)^\alpha-B^\alpha + iu\alpha B^{\alpha-1} \right)}}\nonumber.
\end{align}
The trajectories $t\rightarrow X_t$ can be sampled by simulating the increments on a grid of observation rate $\Delta>0$. The following steps are outlined for $Y^+$ but can be straightforwardly adapted to $Y^{-}$, replacing $P$ and $A$ by $Q$ and $B$.

\subsubsection{Bauemer-Merchaer algorithm}
The Bauemer-Merchaer algorithm is based on the density equation \eqref{eq:density_CTS_stable} that provides a proposal distribution for the acceptance-rejection algorithm (see \cite{kawai_simulation_2011} for more details). However, it applies differently to the finite variation case $\alpha \in (0,1)$ to the infinite variation case $\alpha \in (1,2)$. In the first case, the sampling is exact while in the other it only provides an approximation.

\subsubsection*{Finite variation case, $\alpha \in (0,1)$.}

Given the explicit expression of $\phi_{Y_\Delta^+}$ in \eqref{eq:tempered_subordinator_cf}, the process $Y^+$ is a drifted tempered subordinator, and it can be decomposed into $Y_\Delta^+ = Z_\Delta^+ - \Delta \alpha \Gamma(-\alpha)PA^{\alpha-1} \Delta$. The density $f_{Z_\Delta^+}$ is then merely given after shifting in  \eqref{eq:density_CTS_stable} for $x\in \R$ by
\begin{align}\label{eq:density_tempered_stable_link}
f_{Z_\Delta^+}(x) = e^{-Ax - \Delta P\Gamma(-\alpha)A^\alpha} f_{S_\Delta^+}(x).
\end{align}
Leveraging the fact that the support of $S_\Delta^+$ is $supp(f_{S_\Delta^+})= \R_+$ (see \eqref{eq:support_stable}), the ratio $f_{Z_\Delta^+}(x)/f_{S_\Delta^+}(x)$ can uniformly be bounded by $M= e^{-\Delta P\Gamma(-\alpha)A^\alpha}$. Consequently, $Z_\Delta^+$ can be sampled by accepting  $S_\Delta^+$ as a draw when $U\leq e^{-AS_\Delta^+}$ for $U\sim \mathcal{U}(0,1)$. The increments $S_\Delta^+$ are straightforwardly sampled using Algorithm \ref{algo:simulation_S(sigma,beta,delta)} of Section \ref{sec:simulation_stable}. Algorithm \ref{algo:simulation_Y_t^+} summarizes the above discussion and provides an exact algorithm for simulating increments $Y_\Delta^+$ when $\alpha \in (0,1)$. In Figure \ref{fig:CTS_trajectories_histograms} trajectories for $\alpha=0.5, P=1.7$ and $Q=0.3$ are displayed, with $n=1000$ increments across different time scales ($\Delta \in \{0.01,0.1,1\}$) and computation time is fast. 

\quad

\textbf{Algorithm analysis:} Let $N_+$ be the number of rejected samples in Algorithm \ref{algo:simulation_Y_t^+}. It is known that 
$N^+$ has a geometric distribution of 'success' parameter $s^+$ given by the formula $s^+:=s^+_{\Delta,\alpha,P,A} = \Pro(U \leq e^{-AS_\Delta^+})$. Conditioning by $S_\Delta^+$, yields an expression for $s^+$  
\begin{align}\label{eq:success_CTS_algo_FV}
  s^+= \E\left(e^{-AS_\Delta^+} \right)= e^{\Gamma(-\alpha)\Delta P A^\alpha}.  
\end{align}
The expected number of iterations is  $\E(N^+) = e^{-\Gamma(-\alpha) P\Delta A^{\alpha}}$ ($\Gamma(-\alpha)<0$ when $\alpha \in (0,1)$). Consequently, the computing time grows exponentially with the sampling rate $\Delta>0$ and the jump mass $P$. More precisely, the algorithm is faster in a high-frequency framework; as $s^+$ increases to $1$ and $N^+$ decreases to $0$ when $\Delta \rightarrow 0$. Whereas, when $\Delta \rightarrow +\infty$ the success parameters $s^+$ decrease to $0$ and the number of iterations explodes toward infinity exponentially. As a function of $\alpha$, it appears from \eqref{eq:success_CTS_algo_FV} that it is more costly to sample when $\alpha$ is close to the edges $0$ or $1$. For the bilateral Lévy measure the expected number of iteration is given by $\E(N^+ + N^-) = e^{-\Gamma(-\alpha)P\Delta A^\alpha} + e^{-\Gamma(-\alpha)Q\Delta B^\alpha},$ where $N^-$ is defined analogously to $N^+$ for $Y_\Delta^{-}$.
\subsubsection*{Infinite Variation Case, $\alpha \in (1,2)$}

 In this case, the latter acceptance-rejection method cannot be directly implemented for $\alpha \in (1,2)$. In fact the decomposition into $Z_\Delta^+ - \Delta\alpha\Gamma(-\alpha)PA^{\alpha-1}$ is no longer valid. Moreover, the density equation \eqref{eq:density_CTS_stable} cannot be leveraged to provide another rejection algorithm since the support of $S_\Delta^+$ is now the whole line $\R$ (see \eqref{eq:moment_stable}). However, it can be adapted by injecting $c>0$, which serves as a truncation parameter, into the acceptance condition.  The outcome of such an algorithm is an approximation of the desired distribution $Y_\Delta^+$.
 
In practice a uniform random variable is sampled $U\sim \mathcal{U}([0,1])$ and $Y\sim S_\Delta^+$ then $Y$ is accepted when  $U\leq e^{-A(S_\Delta^+ +  c)}$. The result of this procedure, $\tilde{Y}_{c,\Delta}^+$, is a good candidate for approximating $Y_\Delta^+$. Theoretical guarantees are provided in Theorem 8 of \cite{Baeumer2010TemperedSL} stating that $\tilde{Y}_{c,\Delta}^+ \stackrel{L^1}\rightarrow Y_\Delta^+$ when $c\rightarrow +\infty$. 

\quad
 
\begin{figure}
\begin{center}
$\Delta=0.1$, $\alpha=0.5$ (top), $\alpha=1.5$ (bottom)
\begin{center}
\includegraphics[scale=0.5]{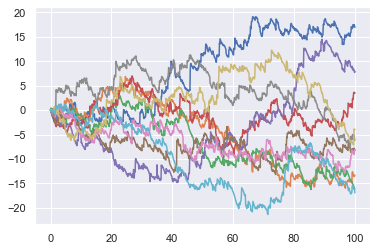}
\includegraphics[scale=0.5]{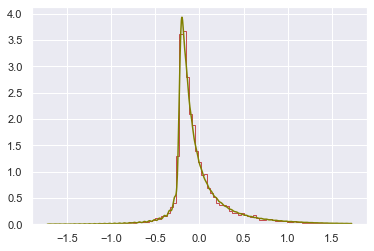}
\includegraphics[scale=0.5]{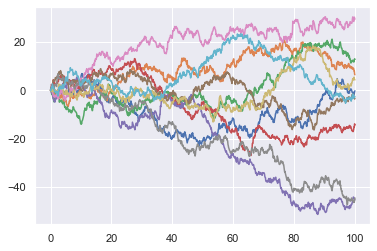}
\includegraphics[scale=0.5]{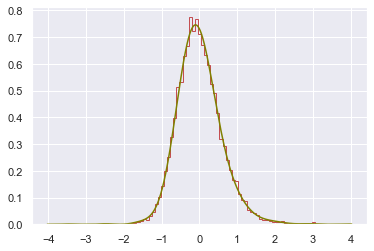}
\end{center}
$\Delta=1$, $\alpha=0.5$ (top), $\alpha=1.5$ (bottom) \begin{center}
\includegraphics[scale=0.5]{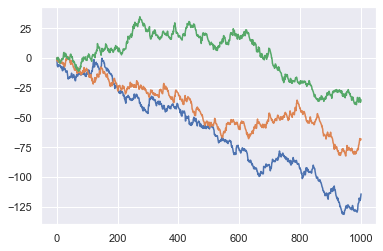}
\includegraphics[scale=0.5]{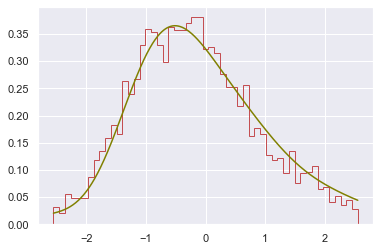} 
\includegraphics[scale=0.5]{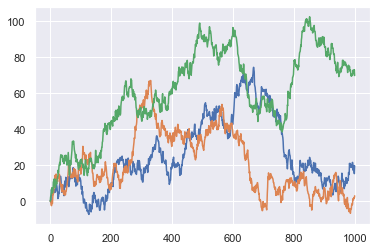}
\includegraphics[scale=0.5]{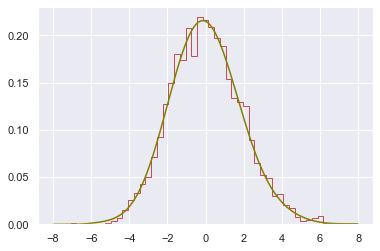}
\end{center}

\caption{Trajectories of a CTS Lévy process $(X_t)_t$ with Lévy triplet $(b_\nu,0,\nu)$, $P=1.7,Q=0.3$, sampled with $n=1000$ observations of step $\Delta \in \{0.01,0.1,1\}$. The parameter $c=1$ for $\Delta<1$ and $c=10$ for $\Delta=1$. Histograms of the increments $X_\Delta$ (red) and numerical approximation CTS density $f_{X_\Delta}$ (green) are displayed. }
\label{fig:CTS_trajectories_histograms}

\end{center}
\end{figure}

\textbf{Algorithm analysis}:
The acceptance rate and the approximation error are interlaced; as $c$ grows, the approximation is better, but the sampling is more costly. The 'success' rate is given by $$s^+=s_{\Delta,\alpha,P,A,c}^+= \Pro(U\leq e^{-A(S_\Delta^+ + c)}) = \E(e^{-A(S_\Delta^+ + c)} \ind_{S_\Delta^+>-c}) + \Pro(S_\Delta^+ \leq -c).$$
The expected number of iterations is then merely the inverse $\E(N^+) = (s^+)^{-1}$. It is clear that $s^+ {\rightarrow} 0$ when $c$ goes to $+\infty$. This suggests that a compromise has to be made in the selection of $c$.
Seen as a function of the step $\Delta>0$, the acceptance rate  $s^+ \underset{\Delta\rightarrow 0}{\rightarrow} 1$, meaning that the high-frequency regime necessitates fewer iterations. Neither the asymptotic behavior of the acceptance rate nor the choice of a suitable approximation metric is straightforward. In \cite{kawai_simulation_2011}, the interplay between the approximation error in Kolmogorov-Smirnov distance $D_{KS}(\tilde{Y}_\Delta^+, Y_\Delta^+)$ and the acceptance rate $s^+$ is numerically explored.  The use of minimization algorithms to numerically compute the parameter $c$ is suggested. In our illustrative examples, the parameter $c>0$ has been selected after preliminary computations on each set of parameters and gives satisfactory results, especially for small $\Delta$. In Figure \ref{fig:CTS_trajectories_histograms} trajectories for $\alpha=0.5, P=1.7$ and $Q=0.3$ are displayed, with $n=1000$ increments across different time scales ($\Delta \in \{0.01,0.1,1\}$). The selected parameters in these cases are $c=1$ for $\Delta\in\{0.01, 0.1\}$ and $c=10$ for $\Delta=1$. The computation time drastically grows with $\Delta$.

\begin{center}
\begin{algorithm}[H]
\textbf{Step 1. }Generate $U\sim \mathcal{U}(0,1)$ and $S_\Delta^{+}\sim S_\alpha(\Delta^{1/\alpha}\sigma,1,0)$.

\textbf{Step 2. } Fix $\tilde c= c \ind_{\alpha \in (1,2)}$ where $c>0$. 

\textbf{Step 3: }\If{$U\leq e^{-A(S_\Delta^+ + \tilde c)}$}
{
\KwRet{ $S_\Delta^+ -\Gamma(1-\alpha)\Delta P A^{\alpha-1}$.}
}
\Else
{Return to Step 1.
}
\caption{Sampling $Y_\Delta^{+}$ for $\alpha \in (0,1)\cup(1,2)$. }
\label{algo:simulation_Y_t^+}
\end{algorithm}
\begin{algorithm}[H]
\textbf{Step 1. } Sample $n$ copies of $(Y_\Delta^+,Y_\Delta^{-})$ using Algorithm \ref{algo:simulation_Y_t^+}.

\textbf{Step 2. } Compute $n$ copies of $X_\Delta = Y_\Delta^{+} - Y_\Delta^{-}$.

\KwRet{$\left(\sum_{j=1}^i X_\Delta^{(i)}\right)_{i=1}^{n}$.}

\caption{Sampling a trajectory of a CTS Lévy process on a grid $\{\Delta,2\Delta,\cdots,n\Delta\}$. }
\label{algo:simulation_Tempered_skeleton)}
\end{algorithm}
\end{center}

\begin{figure}[H]
\begin{center}
$\Delta=0.01$, $\alpha=0.5$ (top), $\alpha=1.5$ (bottom)
\begin{center}
\includegraphics[scale=0.5]{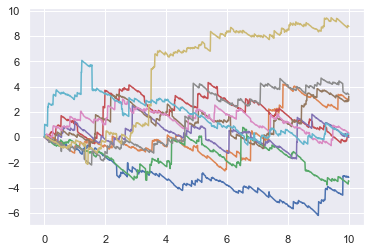}
\includegraphics[scale=0.5]{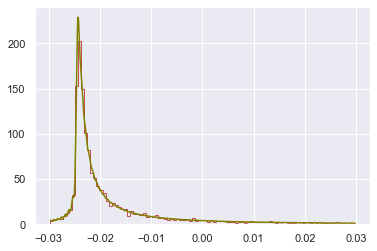}
\includegraphics[scale=0.5]{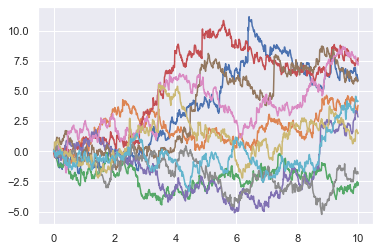}
\includegraphics[scale=0.5]{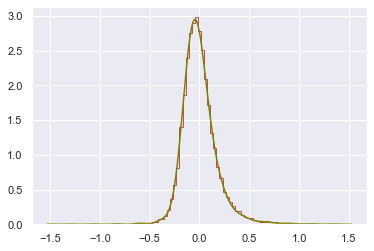}
\end{center}
\end{center}
\end{figure}

\newpage
%%-----------------------------
\nocite{Waskom2021}
\bibliographystyle{abbrv}
\bibliography{main_review_v1}
%%-----------------------------
\end{document}